\documentclass{amsart}
\title{Extensions of generic measure-preserving actions}
\author{Julien Melleray}

\address{Universit\'e de Lyon \\
CNRS UMR 5208 \\
Universit\'e Lyon 1 \\
Institut Camille Jordan \\
43 blvd. du 11 novembre 1918\\
F-69622 Villeurbanne Cedex \\
France}
\email{melleray@math.univ-lyon1.fr}


\usepackage{amssymb}
\usepackage[initials,shortalphabetic]{amsrefs}
\usepackage[all]{xy}

\usepackage{mathpazo}

\usepackage{my-macros}
\numberwithin{equation}{section}

\DeclareMathOperator{\Res}{Res}

\begin{document}

\begin{abstract}
We show that, whenever $\Gamma$ is a countable abelian group and $\Delta \le \Gamma$ is a finitely-generated subgroup, a generic measure-preserving action of $\Delta$ on a standard atomless probability space $(X,\mu)$ extends to a free measure-preserving action of $\Gamma$ on $(X,\mu)$. This extends a result of Ageev, corresponding to the case when $\Delta$ is infinite cyclic.
\end{abstract}
\maketitle
\section{Introduction}
A classical subject in ergodic theory is the study of measure-preserving actions of countable groups on a standard atomless probability space $(X,\mu)$. Given two countable groups $\Delta \le \Gamma$, one may ask whether any measure-preserving action of $\Delta$ on $(X,\mu)$ extends to a meaure-preserving action of $\Gamma$, and it is well-known that, in this generality, the question has a negative answer.

Let $G$ denote the automorphism group of $(X,\mu)$, endowed with its usual Polish group topology. The set of measure-preserving actions of $\Gamma$ on $(X,\mu)$, which is just the set of homomorphisms from $\Gamma$ into $G$, may naturally be identified with a closed subset of $G^{\Gamma}$ (endowed with the product topology), and so the space of $\Gamma$-actions is a Polish space in its own right. Then, one may wonder what happens for a \emph{generic} action of $\Delta$ on $(X,\mu)$, in the sense of Baire category: can a generic action of $\Delta$ on $(X,\mu)$ be extended to an action of $\Gamma$ on $(X,\mu)$? This problem is for instance mentioned in \cite{Kechris2010}*{p. 75}, where one can find an example of countable groups $\Delta \le \Gamma$ such that $\Delta$ is infinite cyclic and a generic action of $\Delta$ cannot be extended to an action of $\Gamma$. 

It is well-known that, for any countable group $\Delta$, generic actions of $\Delta$ are free and so, given a pair of countable groups $\Delta \le \Gamma$, one may also ask whether a generic action of $\Delta$ can be extended to a \emph{free} action of $\Gamma$.
A complete answer to that question, in the case $\Delta=\Z$ and $\Gamma$ is abelian, has been provided by Ageev ~\cite{Ageev2003}*{Theorem 2}: in that case, a generic action of $\Delta$ does extend to a free action of $\Gamma$. In this paper, we extend Ageev's result to the case when $\Delta$ is finitely generated abelian \footnote{Ageev did not publish his proof, so it was unknown to me when writing this article whether his argument was similar to what is presented here. Since then (private commmunication) he told me that his proof was quite different.}, proving the following theorem (Corollary~\ref{c:free} below).

\begin{theorem*}
Let $\Gamma$ be a countable abelian group and $\Delta$ a finitely-generated subgroup of $\Gamma$. Then a generic measure-preserving action of $\Delta$ on a standard atomless probability space $(X,\mu)$ can be extended to a free measure-preserving action of $\Gamma$ on $(X,\mu)$.
\end{theorem*}

Our approach to this question is via \emph{category-preserving maps} and a generalization of the classical Kuratowski--Ulam theorem valid for these maps; these notions were first considered in \cite{Melleray2011}\footnote{After completing a first draft of this paper, I became aware that a similar approach was used by Tikhonov ~\cite{Tikhonov2006} to study embeddings of generic actions of $\Z^d$ in continuous actions of $\R^d$.}.

The paper is organized as follows: first we quickly go over some background on the space of actions of countable groups and properties of generic measure-preserving $\Z$-actions that will be needed in our proof. Next, we recall the definition of a category-preserving map, establish some properties, and discuss the relationship between our approach to the problem tackled here and a classical approach to similar problems in ergodic theory, which is via ``Dougherty's lemma'' and the notion of points of local density for a continuous map between two Polish spaces. Then we give a proof of our main result and discuss possible generalizations. While we prove that most of these possible generalizations are false, we leave open the question of whether one might drop the assumption that $\Delta$ is finitely generated in the statement of our main result.

\emph{Acknowledgements.} When preparing this paper, I benefitted from conversations with several people; in particular I'd like to thank Oleg Ageev for mentioning his paper \cite{Ageev89}; Damien Gaboriau for interesting discussions; Bruno S\'evennec for pointing out the existence of the countable group which is used to prove Theorem~\ref{t:cex}; S\l awomir Solecki for important bibliographical information; and Todor Tsankov for useful remarks and corrections on a first draft of the paper. Work on this project was partially supported by the ANR network AGORA, NT09-461407 and ANR project GRUPOLOCO, ANR-11-JS01-008.

\section{Background and terminology} 
\subsection{The space of actions}
For information on Polish groups and spaces, we refer the reader to \cite{Kechris1995} and \cite{Gao2009a}. The book \cite{Kechris2010} is a good reference for ergodic theory seen from the descriptive set theoretic point of view.

Let $\Gamma$ be a countable group and $G$ be a Polish group. We denote by $\Hom(\Gamma,G)$ the set of homomorphisms of $\Gamma$ into $G$. If we endow $G^{\Gamma}$ with its product Polish topology $\tilde \tau$ (the product of countably many copies of $(G,\tau)$) then $\Hom(\Gamma,G)$ is a closed subset of $G^{\Gamma}$, hence $(\Hom(\Gamma,G), \tilde \tau)$ is a Polish topological space in its own right, and the conjugation action given by $(g \cdot \pi)(\gamma)=g \pi(\gamma)g^{-1}$ is continuous. We may then use Baire category notions in $\Hom(\Gamma,G)$; below, we say that a subset $\Omega \subseteq \Hom(\Gamma,G)$ is \emph{generic}, or \emph{comeager}, if it contains a countable intersection of dense open subsets of $\Hom(\Gamma,G)$. Dually, a set is meager if its complement is comeager; we will often use the formulation ``a generic $\pi \in \Hom(\Gamma,G) $ has property (P)'' to mean ``the set of all $\pi \in \Hom(\Gamma,G)$ which have property (P) is generic''.

By a \emph{standard atomless probability space}, we mean, as usual, a probability space isomorphic to $[0,1]$ endowed with the Lebesgue measure.
If $G=\Aut(X,\mu)$ denotes the automorphism group of a standard atomless probability space $(X,\mu)$, then we endow $G$ with its usual Polish group topology, which is most easily described by its convergent sequences: a sequence $(g_n)$ of elements of $G$ converges to $g \in G$, if, for any measurable $A \subseteq X$, one has $\mu(g_nA \Delta gA) \to 0$ as $n \to + \infty$. Note that two measure-preserving automorphisms of $(X,\mu)$ are identified if they coincide on a set of full measure; since all the subgroups of $G$ that we consider here are countable, this should not cause any confusion and so we'll just neglect sets of measure zero in what follows.

\begin{nota*}
In the remainder of this article, $(X,\mu)$ is a standard atomless probability space and $G$ denotes its automorphism group, endowed with its usual Polish topology. 
\end{nota*}

We say that $\pi \in \Hom(\Gamma,G)$ is \emph{free} if $\mu(\set{x \colon \pi(\gamma)x=x})=0$ for all $\gamma \ne 1$. The set of free measure-preserving actions is dense $G_{\delta}$ in $\Hom(\Gamma,G)$ ~\cite{Glasner1998}. If $\Gamma$ is finite, any two free measure-preserving actions of $\Gamma$ are conjugate, hence a free measure-preserving action of $\Gamma$ on $(X,\mu)$ has a comeager conjugacy class. While such a strong fact is not true in general (for instance, conjugacy classes in $\Hom(\Gamma,G)$ are meager whenever $\Gamma$ is amenable and infinite \cite{Foreman2004}), it is well-known that for any countable group $\Gamma$ there exist actions of $\Gamma$ on $(X,\mu)$ which have a dense orbit under conjugacy ~\cite{Glasner2006a}, so the $0-1$ topological law (see e.g. \cite{Kechris1995}*{8.46}) implies that any conjugacy-invariant subset with the property of Baire is either meager or comeager. This fact seems to be often called the \emph{dynamical alternative}, after the terminology of \cite{Glasner1998}, where it was first established.

Let us also note here a few results that will be useful below; 

\begin{theorem}[King~\cite{King2000}] \label{t:roots}
A generic element of $G$ admits roots of all orders. Actually, the map $g \mapsto g^n$ is category-preserving (see below for a definition) for all $n \ne 0$.
\end{theorem}

This was first proved by King in \cite{King2000}, then de la Rue and de Sam Lazaro gave a simpler presentation of King's proof in \cite{delarue2003}, and improved the above result by showing that a generic element $g$ of $G$ embeds in a flow, i.e there exists a continuous homomorphism $F \colon (\R,+) \to G$ such that $g=F(1)$.

For $g \in G$, we denote by $\langle g \rangle$ the subgroup of $G$ generated by $g$.

\begin{theorem}[Chacon--Schwartzbauer \cite{Chacon1969}] \label{t:centralizer}
For a generic element $g$ of $G$, the centralizer $C(g)$ of $g$ coincides with $\overline{\langle g \rangle}$; in particular $C(g)$ is a maximal abelian subgroup of $G$. 
\end{theorem}
In an earlier version of this article, the above result was incorrectly attributed to King \cite{King1986}; actually, it was proved much earlier: it is stated in \cite{Akcoglu1970}, where the authors say it was already proved by Chacon--Schwartzbauer \cite{Chacon1969} (though the result does not seem to appear explicitly there) \footnote{I am grateful to S. Solecki for pointing this out to me.}. Yet another proof recently appeared in \cite{Melleray2011}.

Stepin and Eremenko, using techniques originated by Ageev and de la Rue--de Sam Lazaro, proved that the centralizer of a generic element is large in the following sense.

\begin{theorem}[Stepin--Eremenko~\cite{Stepin2004}] \label{t:torus}
The infinite-dimensional torus $ \T^{\omega}$ embeds isomorphically (as an abstract group) in the centralizer of a generic element of $G$.
\end{theorem}

What we really need for our proof is a weaker corollary of this theorem, originally proved by Ageev \cite{Ageev2000}, namely the fact that if $g$ is a generic element of $G$ then any finite abelian group isomorphically embeds in the centralizer of $g$.

Though it will not be needed in this paper, we mention for completeness that a recent result of Solecki \cite{Solecki2012}, who proved that the centralizer of a generic element of $\Aut(X,\mu)$  is a continuous homomorphic image of a closed subspace of $L^0(\R)$ and contains an
increasing sequence of finite dimensional tori whose union is dense (in \cite{Solecki2012}, Solecki explains how to use this result to derive the theorem of Stepin--Eremenko quoted above, which is not explicitly proved in \cite{Stepin2004} even though it is stated in the abstract of that paper and can also be derived from the authors' arguments).
\subsection{Category-preserving maps}

As explained in the introduction, our aim is to show that, whenever $\Delta \le \Gamma$ are countable abelian groups and $\Delta$ is finitely generated, a generic element of $\Hom(\Delta,G)$ may be extended to a free element of $\Hom(\Gamma,G)$. In particular, denoting by $\Res$ the restriction map from $\Hom(\Gamma,G)$ to $\Hom(\Delta,G)$, we would like to show that the image of $\Res$ is comeager. Since the conjugacy action of $G$ on $\Hom(\Delta,G)$ is topologically transitive, the $0-1$ topological law implies that it is enough to show that the image of $\Res$ is not meager.

A common approach to this type of question in ergodic theory is based on an observation sometimes called Dougherty's lemma (see e.g \cite{King1986}, \cite{King2000}, \cite{Stepin2004}, \cite{Tikhonov2006}); below we quickly discuss this approach, as well as the technique used in \cite{Melleray2011}, and compare the two.

\begin{defn}
Let $Y,Z$ be topological spaces and $f \colon Y \to Z$ a continuous map. Say that $y \in Y$ is \emph{locally dense} for $f$ if for any neighborhood $U$ of $y$ the set $\overline{f(U)}$ is a neighorhood of $f(y)$.
\end{defn}

\begin{prop}[``Dougherty's lemma '' ]
Assume $Y,Z$ are complete metric spaces, $f \colon Y \to Z$ is continuous and the set of points which are locally dense for $f$ is dense in $Y$. Then $f(Y)$ is not meager.
\end{prop}

\begin{proof}
Assume that the set of points which are locally dense for $f$ is dense and $f(Y) \subseteq \cup_n F_n$ where each $F_n$ is a closed subset of $Z$ with empty interior. Then $Y= \cup_n f^{-1}(F_n)$ so some $f^{-1}(F_n)$ must have nonempty interior, hence must contain a point of local density, so the interior of $\overline{f(f^{-1}(F_n))} \subseteq F_n$ is nonempty, a contradiction.
\end{proof}

Conversely, the next proposition shows that, when $Y$ is separable, the existence of points of local density is necessary for $f(Y)$ to be non meager, though it is certainly not necessary that the set of points of local density be dense; as a side remark, note that this is always a $G_{\delta}$ subset of $Y$, see \cite{King1986}.

\begin{prop}\label{p:nolocaldense}
Assume $Y,Z$ are Polish spaces, $f \colon Y \to Z$ is continuous and let $A=\set{y \in Y \colon y \text{ is not locally dense for } f}$. Then $f(A)$ is meager.
\end{prop}

\begin{proof}
If $y$ is not locally dense for $f$, then there exists an open subset $U$ such that $y \in U$ and $f(y) \not \in \text{Int}(\overline{f(U)})$. 
Hence $f(y) \in \overline{f(U)} \setminus \text{Int}(\overline{f(U)})$. Choosing a countable basis of open subsets $(U_n)$ for the topology of $Y$, we see that
\[
f(A) \subseteq \bigcup_n \overline{f(U_n)} \setminus \text{Int}(\overline{f(U_n)}) \ .
\]
Hence $f(A)$ is meager.
\end{proof}

A different approach, at least on the face of it, was used in \cite{Melleray2011} to study similar questions.

\begin{defn}[\cite{Melleray2011}]  \label{d:catpres}
Let $Y,Z$ be Polish spaces. Say that a continuous map $f \colon Y \to Z$ is \emph{category-preserving} if it satisfies one of the following equivalent conditions:
\begin{enumerate} \romanenum
\item For any comeager $A \subseteq Z$, $f^{-1}(A)$ is comeager.
\item For any nonempty open $U \subseteq Y$, $f(U)$ is not meager.
\item \label{d:catpres3} For any nonempty open $U \subseteq Y$, $f(U)$ is somewhere dense.
\end{enumerate}
\end{defn}

Note that $f$ being category-preserving implies in particular that $f(Y)$ is not meager; also, any continuous open map is category-preserving.

It turns out that the two approaches are equivalent.

\begin{prop}
Assume $Y,Z$ are Polish spaces and $f \colon Y \to Z$ is continuous. Then $f$ is category-preserving if and only if the set of points which are locally dense for $f$ is dense in $Y$.
\end{prop}

\begin{proof}
The implication from right to left is immediate from the definition of a point of local density  and condition~\ref{d:catpres} \eqref{d:catpres3}. To see the converse, assume that there is a nonempty open subset $U$ of $Y$ such that $U$ does not contain any point of local density. Then Proposition~\ref{p:nolocaldense} implies that $f(U)$ is meager, so $f$ is not category-preserving.
\end{proof}

So far, we have explained an approach to showing that the image of the restriction map $\Res \colon \Hom(\Gamma,G) \to \Hom(\Delta,G)$ is comeager; actually, we want to prove more, since we want to prove that a generic action of $\Delta$ extends to a \emph{free} action of $\Gamma$. In other words, we want to prove that the restrictions of free actions of $\Gamma$ form a comeager set in $\Hom(\Delta,G)$. This will come for free (no pun intended) if we etablish that the restriction map is category-preserving: indeed it is easy to check, assuming that the restriction map is category-preserving, that the restrictions of elements taken in any comeager subset of $\Hom(\Gamma,G)$ must form a comeager subset of $\Hom(\Delta,G)$

Saying that $f \colon Y \to Z$ is category-preserving is a strong condition, much stronger than just saying that $f(Y)$ is comeager. This is witnessed by the following theorem.

\begin{theorem}[\cite{Melleray2011}]\label{t:KU}
Let $Y,Z$ be Polish spaces, and $f \colon Y \to Z$ be a category-preserving map. Let also $A$ be a subset of $Y$ with the property of Baire. Then the following assertions are equivalent:
\begin{enumerate} \romanenum
\item $A$ is comeager in $Y$.
\item $\set{z \colon A \cap f^{-1}(z) \text{ is comeager in } f^{-1}(z)}$ is comeager in $Z$.
\end{enumerate}
\end{theorem}

In what follows, we will use the notation $\forall^* y \in Y \ A(y)$ to signify that $A$ is comeager in $Y$; for example, the equivalence in the above theorem, when written using this notation, becomes (for $A \subseteq Y$ with the property of Baire):
\[
(\forall^* y \in Y \ A(y)) \Leftrightarrow (\forall^* z \in Z \ \forall^* y \in f^{-1}(z) \ A(y)) \ .
\]

The above statement, in the case $Y=Y_1 \times Y_2$ and $f$ is a coordinate projection, is the classical Kuratowski--Ulam theorem. It enables one to ``split category along the fibers of $f$'' and will be extremely useful in our proof.

After a first draft of this paper was completed, I became aware of the paper \cite{Tikhonov2006}, where the author considers maps which ``respect genericity'' - i.e such that the inverse image of a comeager set is comeager, and the image of a comeager set is comeager. Any map which respects genericity is category-preserving, and a category-preserving map respects genericity if and only if it has a dense image. Thus, in the context of this article, the two notions are essentially equivalent; compared to \cite{Tikhonov2006}, our approach is somewhat simpler because we have Theorem~\ref{t:KU} at our disposal, while \cite{Tikhonov2006} only uses the classical Kuratowski--Ulam theorem.

The following simple proposition from \cite{Melleray2011} will also be useful. 
 
\begin{prop}\label{p:minimal}
Let $H$ be a Polish group, $Y,Z$ be two Polish $H$-spaces and $f \colon Y \to Z$ a $H$-map. Assume that $Y$ is minimal (i.e, every orbit is dense) and $f(Y)$ is not meager. Then $f$ is category-preserving.
\end{prop}

\begin{proof}
Since $f(Y)$ is not meager, there must exist a point of local density $y \in Y$. Then every point of $H \cdot y$ is a point of local density, so the set of points of local density is dense, and we are done.
\end{proof}

As an example, an immediate corollary of this lemma is that, whenever $H,K$ are Polish groups, a continuous homomorphism $\phi \colon H \to K$ with non-meager image must be category-preserving; actually it is well-known that under these assumptions $\phi$ must be open. Let us give a slightly more interesting example.

\begin{lemma}\label{l:finitegroup}
Let $\Delta \le \Gamma$ be two groups, and assume that $\Delta$ is finite and $\Gamma$ is countable. Then the restriction map $\Res \colon \Hom(\Gamma,G) \to \Hom(\Delta,G)$ is category-preserving.
\end{lemma}

\begin{proof}
Let $Y$ denote the set of all free actions of $\Gamma$ with a dense conjugacy class. $\Res(Y)$ is conjugacy-invariant and contains a free action of $\Delta$, hence it must contain all of them (recall that they are all conjugate) and so is comeager. The action of $G$ on $Y$ is minimal, so $\Res \colon Y \to \Hom(\Delta,G)$ is category-preserving. Since $Y$ is comeager in $\Hom(\Gamma,G)$ (\cite{Glasner2006a}, see \cite{Kechris2010}*{Theorem~10.7} for a simple proof) this implies that $\Res$ is category-preserving.
\end{proof}

In the proof of the main result, we will also need the following facts from \cite{Melleray2011}.

\begin{lemma}[\cite{Melleray2011}]\label{l:monothetic}
Let $H$ be a Polish group and $\Gamma$ be a countable abelian group such that, for a generic $\pi \in \Hom(\Gamma \times \Z,H)$, $\overline{\pi(\Gamma \times \Z)}=\overline{\pi(\Gamma)}$. Then the restriction map $\Res \colon \Hom(\Gamma \times \Z,H) \to \Hom(\Gamma,H)$ is category preserving.
\end{lemma}

This criterion was applied to obtain the following result (also proved, earlier and independently, by Tikhonov \cite{Tikhonov2006}, using his notion of map respecting genericity).

\begin{lemma}[\cite{Melleray2011}]\label{l:products}
For any integers $d \le k$, the restriction map from $\Hom(\Z^k, G)$ to $\Hom(\Z^d,G)$ is category-preserving.
\end{lemma}

\section{Proof of the main result}
We are now ready to give the proof of our main result. We break down the argument in a series of lemmas, most of which are special cases of the main result. Some of the lemmas were already known, but the proofs we give here seem simpler and we try to limit the appeal to black boxes as much as possible; the two results that are used below and that are not proved here or in \cite{Melleray2011} are Theorems~\ref{t:roots} and~\ref{t:torus}. We begin by setting the notation to be used in the proof.

\begin{nota*} We recall that $(X,\mu)$ denotes a standard atomless probability space and $G$ its automorphism group, endowed with its usual Polish topology. All groups are noted multiplicatively (in particular, $1$ stands for the neutral element). If $\Gamma$ is a group and $A$ is a subset of $\Delta$, we denote by $\langle A \rangle$ the group generated by $A$. By $C_k(G)$ we mean the set of commuting $k$-uples of elements of $G$, which we identify with $\Hom(\Z^k, G)$ whenever it is convenient. If $\bar g \in C_k(G)$, $C(\bar g)$ denotes the centralizer of $\langle \bar g \rangle $; when $\Gamma$ is some finite abelian group and $\bar g \in C_k(G)$, we let $C_{\Gamma}(\bar g)$ denote the set of $\Gamma$-actions on $(X,\mu)$ which commute with $\bar g$. Similarly, if $\Gamma$ is some finite group, $\pi \in \Hom(\Gamma,G)$ and $k \in \omega$, $C_k(\pi)$ stands for the set of commuting $k$-uples of elements of $G$ which commute with $\pi(\gamma)$ for all $\gamma \in \Gamma$. Whenever $d \le k$ and we view $\Z^d$ as a subgroup of $\Z^k$, we view it as the subgoup of $k$-uples generated by the first $d$ elements of the natural basis of $\Z^k$.
\end{nota*}

\begin{lemma} \label{proof:lem1}
Let $d$ be an integer, and $\Gamma$ be a finite abelian group. Then the restriction map $\Res \colon \Hom(\Z^d \times \Gamma,G) \to \Hom(\Z^d,G)$ is category-preserving. 
\end{lemma}

\begin{proof}
Since the restriction map $\Hom(\Z^d, G) \to \Hom(\Z,G)$ is category-preserving (Lemma~\ref{l:products}), Theorems~\ref{t:centralizer} and~\ref{t:torus} imply that 
\[
\forall^* \bar g \in C_d(G) \ C(\bar g)= \cl{\langle g_1 \rangle} \text{ is abelian and contains an isomorphic copy of }  \T^{\omega} \ .
\] 
It is also well-known that a generic element of $\Hom(\Z^d, G)$ is free ergodic. Now, choose $\bar g \in C_d(G)$ such that the action $\Z^d \actsl (X,\mu)$ induced by $\bar g$ is free ergodic, and $C(\bar g)= \cl{\langle g_1 \rangle}$ is an abelian group containing an isomorphic copy of $ \T^{\omega}$.

We may extend the action $\Z^d \actsl (X,\mu)$ induced by $\bar g$ as follows: pick an isomorphic embedding $\phi \colon \Gamma \to C(\bar g)$, and set $\pi(n_1,\ldots,n_d,\gamma)=g_1^{n_1}\cdots g_d^{n_d}\phi(\gamma)$. This action extends the original action of $\Z^d$ to $\Z^d \times \Gamma$, and we claim that $\pi$ is free. To see this, pick $h \in \langle \bar g \rangle$ and $\gamma \in \Gamma$, and assume that 
\[ \mu(\set{x \colon h\phi(\gamma)(x)=x})>0 \ .
\] 
Since $\set{x \colon h \phi(\gamma)(x)=x}$ is $\langle \bar g \rangle$-invariant and $\langle \bar g \rangle$ acts ergodically, we get 
\[\mu(\set{x \colon h\phi(\gamma)(x)=x})=1 \ .
\]
Since $\gamma$ has finite order, $h$ must be of finite order and so $h=1$ since the original $\Z^d$-action was free. Hence $\phi(\gamma)=1$, so $\gamma=1$.

Thus a generic $\Z^d$-action extends to a free action of $\Z^d \times \Gamma$; aplying Proposition~\ref{p:minimal} as in the proof of Lemma~\ref{l:finitegroup}, we obtain that $\Res$ is category preserving.
\end{proof}

An immediate corollary of this and Lemma \ref{l:monothetic} is the following.

\begin{lemma}\label{proof:lem2}
Let $\Gamma$ be a finite abelian group, and $d \le k$ be two integers. Then the restriction map $\Res \colon \Hom(\Z^k \times \Gamma,G) \to \Hom(\Z^d \times \Gamma,G)$ is category-preserving.
\end{lemma}

\begin{proof}
It is enough to show the above result when $k=d+1$ (a composition of category-preserving maps is still category-preserving, a fact that will be useful to us more than once) and $d \ge 1$, since the case $d=0$ is covered by Lemma~\ref{l:finitegroup}. From Lemma~\ref{proof:lem1} and Theorem~\ref{t:centralizer} we obtain 
\[
\forall^* \pi \in \Hom(\Z^{d+1} \times\Gamma ,G) \ \cl{\pi(\Z)} \text{ is maximal abelian} \ .
\]
Hence
\[
\forall^* \pi \in \Hom(\Z^{d+1} \times \Gamma,G) \ \cl{\pi(\Z^{d+1} \times \Gamma)}= \cl{\pi(\Z)}=\cl{\pi(Z^{d}\times \Gamma)} \ .
\]
Applying Lemma \ref{l:monothetic} yields the desired result.
\end{proof}

\begin{lemma}\label{proof:lem3}
Fix an integer $d$. Then, for a generic $\bar g \in C_d(G)$, the group $C(\bar g)$ is divisible.
\end{lemma}

\begin{proof}
Fixing some nonzero integer $p$, it is enough to show that for a generic $\bar g$ any element of $C(\bar g)$ has a $p$-th root in $C(\bar g)$. We begin with the case when $d$ is equal to $1$. We know, by Theorem~\ref{t:roots}, that 
\[
\forall^* h \in G \ \exists f \in G \ h=f^p \ .
\]
Applying Lemma~\ref{l:products} and the fact that the centralizer of a generic $g$ coincides with $\cl{\langle g \rangle}$, we obtain
\[
\forall^* (g,h) \in C_2(G) \ C(h)=C(g)=\cl{\langle g \rangle} \text{ and } \exists f \in G \ h=f^p \ .
\]
Since the equation $h=f^p$ implies that $f$ and $h$ commute, the above equation implies, using Lemma~\ref{l:products} again (projecting on the other coordinate) and Theorem~\ref{t:KU}, that
\[
\forall^* g \in G \ \forall^* h \in C(g)=\cl{\langle g \rangle} \ \exists f \in C(g) \ h=f^p \ .
\]
In other words, for a generic $g$ and any integer $p$ the homomorphism $h \mapsto h^p$ of the abelian Polish group $C(g)$ has a comeager image, hence it must be surjective, and we are done.

The general case immediately follows from Lemma~\ref{l:products}, since
\[
\forall^* \bar g \in C_d(G) \ C(\bar g)= C(g_1) \text{ and } C(g_1) \text{ is divisible.}
\]
\end{proof}

The following closely related lemma was already proved in \cite{Tikhonov2006} (using different vocabulary).

\begin{lemma}\label{proof:lem5}
For any nonzero integers $(n_1,\ldots,n_k)$ the map $(g_1,\ldots,g_k) \mapsto (g_1^{n_1},\ldots,g_k^{n_k})$ is category-preserving from $C_k(G)$ to itself.
\end{lemma}

\begin{proof}
It is enough to show that the map $(g_1,\ldots,g_k) \mapsto (g_1,\ldots,g_{k-1},g_k^p)$ is category-preserving for each $k \ge 2$ (the case $k=1$ is the content of Theorem~\ref{t:centralizer}) and $p \ne 0$. Let $O$ be comeager in $C_k(G)$. 
We have 
\[
\forall^* (g_1,\ldots,g_{k-1}) \in C_{k-1}(G) \ \forall^* h \in C(g_1,\ldots,g_{k-1}) \ (g_1,\ldots,g_{k-1},h) \in O \ .
\]
From Lemma~\ref{proof:lem3} we know that $h \mapsto h^p$ is a surjective homomorphism from $C(g_1,\ldots,g_{k-1})$ to itself for a generic $(g_1,\ldots,g_{k-1}) \in C_{k-1}(G)$, so it is category-preserving and the above equation yields
\[
\forall^* (g_1,\ldots,g_{k-1}) \in C_{k-1}(G) \ \forall^* h \in C(g_1,\ldots,g_{k-1})  \ (g_1,\ldots,g_{k-1},h^p) \in O \ .
\]
That is,
\[
\forall^* (g_1,\ldots,g_k) \in C_k(G) \ (g_1,\ldots,g_{k-1},g_k^p) \in O \ .
\]
\end{proof}

\begin{lemma}\label{proof:lem4}
Let $\Gamma$ be a finite abelian group, and $\Delta$ be a torsion-free subgroup of $\Z^k \times \Gamma$. Then the restriction map 
$\Res \colon \Hom(\Z^k \times \Gamma,G) \to \Hom(\langle \Delta,\Gamma \rangle,G)$ is category-preserving.
\end{lemma}

Note before the proof that, as a group, $\langle \Delta,\Gamma \rangle \cong \Delta \times \Gamma$; also, $\Delta \cong \Z^d$ for some $d \le k$.

\begin{proof}
We may assume that $\Delta$ is nontrivial, i.e.~ isomorphic to $\Z^d$ for some $d \in \set{1,\ldots,k}$. Using the existence of simultaneous bases for subgroups of a free finitely-generated abelian group, we may find a basis $e_1,\ldots,e_k$ of $\Z^k$, nonzero integers $n_1,\ldots,n_d$ and elements $\gamma_1, \ldots \gamma_d$ of $\Gamma$ such that 
\[
\Delta= \langle (e_i^{n_i},\gamma_i) \colon i \le d \rangle \ .
\]
Using these particular generating sets, the restriction map may be identified with the map 
\[ \begin{cases}
         \Hom(\Z^k \times \Gamma,G) \to \Hom(\Z^d \times \Gamma,G) &  \\
        (g_1,\ldots,g_k, \pi) \mapsto (g_1^{n_1}\pi(\gamma_1),\ldots,g_d^{n_d}\pi(\gamma_d),\pi) & \\
        \end{cases}      
 \] 

Let $O$ be a comeager subset of $\Hom(\Z^d \times \Gamma,G)$. We know from Lemma~\ref{l:finitegroup} and Theorem~\ref{t:KU} that 
\[
\forall^* \pi \in \Hom(\Gamma,G) \ \forall^* \bar g \in C_d(\pi) \ (\bar g, \pi) \in O \ .
\]
Since, for each fixed  $\pi$, the map $\bar g \mapsto (g_1 \pi(\gamma_1),\ldots,g_d \pi(\gamma_d))$ is a homeomorphism of $C_d(\pi)$, this may be rewritten as 
\[
\forall^* \pi \in \Hom(\Gamma,G) \ \forall^* \bar g \in C_d(\pi) \ (g_1 \pi(\gamma_1),\ldots,g_d\pi(\gamma_d),\pi) \in O \ .
\]
Using Theorem~\ref{t:KU} and Lemma~\ref{proof:lem1}, this is the same as 
\[
\forall^* \bar g \in C_d(G) \ \forall^* \pi \in C_{\Gamma}(\bar g) \ (g_1 \pi(\gamma_1),\ldots,g_d\pi(\gamma_d),\pi) \in O \ .
\]
Applying Lemma~\ref{proof:lem5} this yields
\[
\forall^* \bar g \in C_d(G) \ \forall^* \pi \in C_{\Gamma}(\bar g) \ (g_1^{n_1} \pi(\gamma_1), \ldots,g_d^{n_d} \pi(\gamma_d),\pi) \in O  \ .
\]
Finally, Lemma~\ref{proof:lem2} and Theorem~\ref{t:KU} lead to the desired
\[
\forall^* \pi \in \Hom(\Z^k \times \Gamma,G) \ \Res(\pi) \in O \ .
\]
\end{proof}

\begin{lemma}\label{proof:lem6}
Let $\Gamma_2 \le \Gamma_1$ be finite abelian groups and $d$ be an integer. Then $\Res \colon \Hom(\Z^d \times \Gamma_1,G) \to \Hom(\Z^d \times \Gamma_2,G)$ is category-preserving.
\end{lemma}

\begin{proof}
We claim that, for a generic $\bar g \in C_d(G)$, any homomorphism from $\Gamma_2$ to $C(\bar g)$ extends to a homomorphism from $\Gamma_1$ to $ C(\bar g)$ . 
The proof of this is the same as the classical proof of the fact that a character of a finite abelian group extends to a character of a finite abelian supergroup. By induction, it is enough to prove the above fact in case there exists $x \in \Gamma_1$ such that $\Gamma_1= \langle \Gamma_2,x \rangle$.
 Pick such an $x$, and let $m$ be the smallest nonnegative integer such that $x^m=\gamma \in \Gamma_2$. Let also $\bar g$ be such that $C(\bar g)=\cl{\langle \bar g \rangle}$ is abelian and divisible.
 
Then for any homomorphism $\phi$ from $\Gamma_2$ to $C(\bar g)$, we may find $h \in C(\bar g)$ such that $h^m=\phi(\gamma)$, and extend $\phi$ to $\Gamma_1$ by setting, for any $k < \omega$ and $\delta \in \Gamma_2$, $\phi(x^k \delta)=h^k \phi(\delta)$. The verification that this is a well-defined homomorphism from $\Gamma_1$ to $C(\bar g)$ is left to the reader.

Since $\Gamma_1, \Gamma_2$ and $C(\bar g)$ are abelian, the restriction map from $C_{\Gamma_1}(\bar g)$ to $C_{\Gamma_2}(\bar g)$ is a (continuous) homomorphism between two Polish groups, and we just proved that this homomorphism is surjective: hence it must be open, and in particular category-preserving.

We are almost done: let $O$ be a comeager subset of $\Hom(\Z^d \times \Gamma_2,G)$. Then we know from Lemma~\ref{proof:lem2} that 
\[
\forall^* \bar g  \in C_d(G)\ \forall^* \pi \in C_{\Gamma_2}(\bar g) \ (\bar g, \pi) \in O \ .
\]
Since a generic $\bar g \in C_d(g)$ is such that $C(\bar g)= \cl{\langle \bar g \rangle}$ is abelian and divisible, our reasoning above yields
\[
\forall^* \bar g \in C_d(G) \ \forall^* \phi \in C_{\Gamma_1}(\bar g) \ (\bar g, \phi_{|\Gamma_2}) \in O \ .
\]
That is, 
\[
\forall^* \pi \in \Hom(\Z^d \times \Gamma_1,G) \ \Res(\pi) \in O \ .  
\]

\end{proof}

Now to the proof of our main result.

\begin{theorem}\label{t:main}
Let $\Gamma$ be a countable abelian group and $\Delta$ be a finitely generated subgroup of $\Gamma$. Then the restriction map $\Res \colon \Hom(\Gamma,G) \to \Hom(\Delta,G)$ is category-preserving.
\end{theorem}

\begin{proof}
We begin with the case when $\Gamma$ is finitely generated.  We may assume that for some nonzero integers $d \le k$, finite groups $F_1 \le F_2$, and a subgroup $H$ of $\Z^k \times F_2$ isomorphic to $\Z^d$, we have 
\[
\Gamma= \Z^k \times F_2 \text{ and } \Delta=H \times F_1 \ .
\]
Then aplying Lemmas~\ref{proof:lem4} and~\ref{proof:lem6} to the following sequence yields the desired result: 
\[\Hom(Z^k \times F_2,G) \to \Hom(\langle H, F_2 \rangle, G) \cong \Hom(H \times F_2,G) \to \Hom(H \times F_1,G) \ .
\]

Now we turn to the case when $\Gamma$ is not finitely generated; we apply the same method as in \cite{Melleray2011} to deal with that case. Going back to the definition of a category-preserving map, we pick a dense open
 $O \subseteq \Hom(\Delta,G)$ and a nonempty open $U \subseteq \Hom(\Gamma,G)$; without loss of generality we may assume that there is a finitely generated subgroup $\Gamma'$ of $\Gamma$ containing $\Delta$, and a nonempty open subset $U'$ of $\Hom(\Gamma',G)$, such that
 \[
 U=\set{\pi \in \Hom(\Gamma,G) \colon \pi_{|\Gamma'} \in U'} \ .
 \]
Since $\Gamma'$ is finitely generated, the restriction map from $\Hom(\Gamma',G)$ to $\Hom(\Delta,G)$ is category-preserving, so 
\[U''= \set{\pi \in U' \colon \pi_{|\Delta} \in O} 
\] is open nonempty in $\Hom(\Gamma',G)$.
Letting $\gamma_1,\ldots,\gamma_n$ denote generators of $\Gamma'$, we may further assume that there exist $\pi_0 \in U''$, $\eps>0$ and a finite measurable partition $\mcA$ of $(X,\mu)$ such that, for any $\pi \in \Hom(\Gamma',G)$, one has
\[
\pi \in U'' \Leftrightarrow \forall A \in \mcA \ \forall i \in \set{1,\ldots,n} \ \mu(\pi_0(\gamma_i)(A)\Delta \pi(\gamma_i)(A)) < \eps \ . 
\]
Let $\psi \colon \Gamma \actsl (X^{\Gamma/\Gamma'}, \mu^{\Gamma/\Gamma'})$ be the action of $\Gamma$ co-induced by $\pi_0$ (see \cite{Kechris2010}*{10(G)}), and let $\Theta \colon X^{\Gamma/\Gamma'} \to X $ be defined by 
$\Theta(f)=f(\Gamma')$. Then, denoting by $\mcB$ the (finite) measurable partition generated by $\mcA$ and $\pi_0(\mcA)$, one may pick a measure-preserving bijection $T \colon X^{\Gamma/\Gamma'} \to X$ such that $T^{-1}(B)=\Theta^{-1}(B)$ for all $B \in \mcB$, and define $\phi \in \Hom(\Gamma,G)$ by setting $\phi(\gamma)=T \psi(\gamma) T^{-1}$. The construction ensures that $\phi_{|\Gamma'} \in U''$, so $\phi \in U \cap \Res^{-1}(O)$, and we are done.
\end{proof}

\begin{remark*}
In general, given a Polish group $H$, and countable groups $\Delta \le \Gamma$, one may ask whether the restriction map from $\Hom(\Gamma,H)$ to $\Hom(\Delta,H)$ is category-preserving. The current paper is concerned with the case $H=\Aut(X,\mu)$, but this question is certainly also of interest for other Polish groups, such as the isometry group of the Urysohn space. Thus is seems worth pointing out that the strategy of proof presented here could conceivably be adapted to other cases; looking at the proof, it is clear that the most important (and, probably, hardest) step is to understand what happens when $\Gamma=\Z$ and $\Delta=n \Z$. So, Chacon--Schwartzbauer's result (Theorem \ref{t:roots}) seems to be the key to understand the general situation.
If one were able to show that a generic isometry of the Urysohn space admits infinitely many $n$-th roots for any $n$ then, using a line of reasoning similar to the one presented here (and some results of \cite{Melleray2011}), one could prove that Theorem \ref{t:main} holds also when $G$ is the isometry group of the Urysohn space.

\end{remark*}

\begin{cor} \label{c:free}
Let $\Gamma$ be a countable abelian group and $\Delta$ be a finitely generated subgroup of $\Gamma$. Then a generic measure-preserving action of $\Delta$ extends to a free measure-preserving action of $\Gamma$.
\end{cor}

\begin{proof}
The restriction map $\Res \colon \Hom(\Gamma,G) \to \Hom(\Delta,G)$ is category-preserving and has a non meager, conjugacy-invariant range, hence a comeager range. Thus the image of any comeager subset of $\Hom(\Gamma,G)$ is comeager in $\Hom(\Delta,G)$; as free actions of $\Gamma$ form a dense $G_{\delta}$ subset of $\Hom(\Gamma,G)$, we are done.

\end{proof}

Now, it is reasonable to wonder to which extent one can strengthen Theorem~\ref{t:main}. 
Clearly, if $\Gamma$ is no longer assumed to be abelian, then the restriction map does not need to be category-preserving in general - indeed, it will never be category-preserving if $\Delta$ is an abelian central subgroup of $\Gamma$, $\Gamma$ is nonabelian and $\Delta$ contains an element with infinite order. To see this, note that it follows from our results that for such a $\Delta$, a generic $\pi \in \Hom(\Delta,G)$ is such that $C(\pi)=\cl{\pi(\Delta)}$ is an abelian group; thus, for a generic $\pi \in \Hom(\Delta,G)$, no action $\tilde \pi$ of $\Gamma$ extending $\pi$ can be free, since $\tilde{\pi}(\Gamma) \le C(\pi)$ must be abelian. Hence, even if $\Delta$ is infinite cyclic and $H$ is a finite nilpotent group, the restriction map $\Res \colon \Hom(\Delta \times H,G) \to \Hom(\Delta,G)$ need not be category-preserving (the simplest counterexample being $\Gamma=\Z \times Q_8$, where $Q_8$ denotes the quaternion group, and $\Delta=\Z$).

Then, one may wonder whether a weakening of Theorem~\ref{t:main} holds in greater generality, namely, whether when $\Delta \le \Gamma$ and $\Delta$ is (say) abelian, a generic measure-preserving action of $\Delta$ may be extended to an action of $\Gamma$ (we are no longer asking that the extension be free). This statement is also false in general; a counterexample appears in \cite{Kechris2010}*{p. 75}. We quickly mention another counterexample: Ageev proved in \cite{Ageev89} that a generic element of $G$ is not conjugate to its inverse, thus a generic measure-preserving $\Z$-action cannot be extended to an action of any nontrivial semidirect product $\Z \rtimes H$. Indeed, the only nontrivial group automorphism of $\Z$ is $n \mapsto -n$, so the image of the generator of $\Z$ under any morphism from $\Z \rtimes H$ to $G$ must be conjugate to its inverse, and generically this does not happen.

In view of the above counterexamples, it is tempting to ask the following question: if $\Gamma$ is a countable group and $\Delta$ is an abelian central subgroup of $\Gamma$, is it true that a generic measure-preserving action of $\Delta$ may be extended to a measure-preserving action of $\Gamma$? This is false even under the additional assumption that $\Delta$ is solvable (I do not know of an example with $\Delta$ nilpotent but it seems likely that such an example exists); I am grateful to Bruno S\'evennec for pointing out to me the existence of a group $\Gamma$ as in the example below.

\begin{theorem}\label{t:cex}
There exists a countable polycyclic (hence, solvable) group $\Gamma$, whose center contains an infinite cyclic subgroup $\Delta$, such that a generic action of $\Delta$ cannot be extended to an action of $\Gamma$. 
\end{theorem}

\begin{proof}
Let $\Gamma$ be the countable group defined by the following presentation with three generators:
\[
\Gamma= \langle x,y,z|\ \left[x,y \right]=z^4, \  zxz^{-1}=x^{-1}, \ z y z^{-1}=y^{-1} \rangle \ .
\]
K.A. Hirsch showed that $\Gamma$ is torsion-free, polycyclic, and its abelianization is finite (see \cite{Hirsch1952}  for the original reference; this is exercise 15 p.158 of \cite{Robinson1996}); from the above presentation we see that $z^2$ is central.

Let $\Delta=\langle z^2 \rangle \le \Gamma$. We know that, for a generic action $\pi$ of $\Delta$, $\pi$ is free and the centralizer of $\pi(z^2)$ is an abelian group; since $z^2$ is central in $\Gamma$, any action $\tilde \pi$ of $\Gamma$ extending $\pi$ must take its values in the abelian group $C(\pi(z^2))$, hence have finite range since $\Gamma$ has no infinite abelian quotients. Thus such a $\tilde \pi$ cannot exist, since $\pi$ is free and so $\tilde \pi(\Delta)=\pi(\Delta)$ must be infinite if $\tilde \pi$ extends $\pi$.
\end{proof}

While it would be interesting to understand exactly for which pairs of groups $\Delta \le \Gamma$ the restriction map from $\Hom(\Gamma,G)$ to $\Hom(\Delta,G)$ is category-preserving, or simply has a comeager image, the above examples show that a general answer will necessarily be complicated and is perhaps too much to ask for. Still, perhaps one may generalize Theorem~\ref{t:main} to the case when $\Delta$ is any countable abelian group (i.e, is no longer supposed to be finitely generated), and we end the paper with that question.

\begin{question}
Is it true that, whenever $\Delta \le \Gamma$ are countable abelian groups, the restriction map $\Res \colon \Hom(\Gamma,G) \to \Hom(\Delta,G)$ is category-preserving?
\end{question}

\bibliography{mybiblio}

\end{document}